\documentclass[12pt,reqno,b5paper]{amsproc}
\usepackage{color,anysize,float}
\usepackage{citeref}
\usepackage{txfonts,esint}
\usepackage[colorlinks,urlcolor=blue,citecolor=blue,linkcolor=blue]{hyperref}

\usepackage{graphicx}%

\marginsize{15mm}{15mm}{1cm}{1cm}
\renewcommand{\bibitempages}[1]{}

\newtheorem{thm}{\sc Theorem}[section]

\newtheorem{lem}[thm]{\sc Lemma}
\newtheorem{prop}[thm]{\sc Proposition}
\theoremstyle{definition}

\theoremstyle{remark}
\newtheorem{rem}[thm]{\it Remark}
\newtheorem*{pf}{\it Proof}

\numberwithin{equation}{section}
\allowdisplaybreaks

\AtEndDocument{\bibliographystyle{siam}
\bibliography{C:/bib/all}}

\begin{document}

\title[Quasilinear Schr\"{o}dinger equations]{Quasilinear Schr\"{o}dinger equations with concave and convex nonlinearities}
\author[S. Liu and L.-F. Yin]{Shibo Liu\,$^\mathrm{a}$, Li-Feng Yin\,$^\mathrm{b}$\vspace{-1em}}
\dedicatory{$^\mathrm{a}$Department of Mathematical Sciences, Florida Institute of Technology\\
Melbourne, FL 32901, USA\\
$^\mathrm{b}$School of Mathematical Sciences, Xiamen University\\
Xiamen 361006, P.R. China}
\thanks{\today\hfill Emails: \tt sliu@fit.edu, yin136823@163.com}
\maketitle
\begin{abstract}
In this paper, we consider the following quasilinear Schr\"{o}dinger equation
\begin{align*}
-\Delta u-u\Delta(u^{2})=k(x)\left\vert u\right\vert ^{q-2}u-h(x)\left\vert
u\right\vert ^{s-2}u\text{,\quad}u\in D^{1,2}(\mathbb{R}^{N})\text{,}
\end{align*}
where $1<q<2<s<+\infty$. Unlike most results in the literature, the exponent $s$ here is allowed to be supercritical $s>2\cdot2^{\ast}$. By taking advantage of geometric properties of a nonlinear transformation $f$ and a variant of Clark's theorem, we get a sequence of solutions with negative energy in a space smaller than $D^{1,2}(\mathbb{R}^{N})$. Nonnegative solution at negative energy level is also obtained.
\end{abstract}

\section{Introduction}

In this paper we consider quasilinear stationary  Schr\"{o}dinger equations of
the form%
\begin{equation}
\left\{
\begin{array}
[c]{l}%
-\Delta u-u\Delta(u^{2})=k(x)\left\vert u\right\vert ^{q-2}u-h(x)\left\vert
u\right\vert ^{s-2}u\text{,}\\
u\in D^{1,2}(\mathbb{R}^{N})\text{,}%
\end{array}
\right.  \label{e9}%
\end{equation}
where $1<q<2<s<\infty$. This kind of equations arise when we are looking for
standing waves $\psi(t,x)=\mathrm{e}^{-i\omega t}u(x)$ for the time dependent
quasilinear Schr\"{o}dinger equation%
\[
i\psi_{t}=-\Delta\psi-\psi\Delta(\left\vert \psi\right\vert ^{2})-\tilde{g}(x,\left\vert
\psi\right\vert ^{2})\psi\text{,\qquad}\left(  t,x\right)  \in\mathbb{R}%
\times\mathbb{R}^{N}\text{.}%
\]
Quasilinear Schr\"{o}dinger equations have captured great interest in the last
two decades because they model several important physical phenomena including
superfluid film in plasma physics, self-trapped electrons in quadratic or
hexagonal lattices, see \cite{MR1899450,MR2630099} and references therein for
more details.

The problem (\ref{e9}) possesses a variational structure. Formally, it is the
Euler--Lagrange equation of the functional%
\[
J(u)=\frac{1}{2}\int\left(  1+2u^{2}\right)  \left\vert \nabla u\right\vert
^{2}-\frac{1}{q}\int k(x)\left\vert u\right\vert ^{q}+\frac{1}{s}\int
h\left\vert u\right\vert ^{s}\text{,}%
\]
where from now on all integrals are taken over $\mathbb{R}^{N}$ unless stated
explicitly. However, $J$ can only be defined on a proper subset of
$D^{1,2}(\mathbb{R}^{N})$, hence the standard variational methods could not be
applied. To overcome this difficulty, Liu \emph{et al.}\ \cite{MR1949452} and
Colin--Jeanjean \cite{MR2029068} introduced a nonlinear transformation $f$
which converts the quasilinear problem into a semilinear one, and enables
us to work with a $C^{1}$-functional
\[
\Phi(v)=J(f(v))=\frac{1}{2}\int\left\vert \nabla v\right\vert ^{2}-\frac{1}%
{q}\int k\left\vert f(v)\right\vert ^{q}+\frac{1}{s}\int h\left\vert
f(v)\right\vert ^{s}
\]
defined on the whole Sobolev space. Since
then, many results about quasilinear Schr\"{o}d\-i\-n\-ger equations appear, mainly
for $4$-superlinear nonlinearities, see \cite{MR2580150,MR3620755,MR3843290,MR3982327}.

In this paper, we study quasilinear Schr\"{o}dinger equations whose
nonlinearity is a combination of concave and convex terms. Elliptic boundary value problems
involving concave and convex terms have attracted great attention since the
pioneering work of Ambrosetti--Brezis--Cerami \cite{MR1276168} and
Bartsch--Willem \cite{MR1301008} on semilinear problems on bounded domain. However, relatively less have been done for quasilinear
Schr\"{o}dinger equations. It seems that \cite{MR2461565} is the first work in
this direction, see \cite{MR4010995} for a more recent result. To apply variational methods, in \cite{MR2461565,MR4010995} and most
papers on quasilinear Schr\"{o}dinger equations, the nonlinearity $g(x,u)$ can at most grow critically, that is%
\[
\vert g(x,u)\vert \le C\big(1+\vert u\vert ^{2\cdot2^{\ast}-1}\big)\text{,}%
\]
here $2^{\ast}=2N/\left(  N-2\right)  $ is the critical
Sobolev exponent. It was pointed out in \cite[Remark 3.13]{MR1949452} that the exponent $2\cdot2^\ast$ behaves like a critical exponent for (\ref{e9}). The nonlinearities in the above mentioned papers on quasilinear Schr\"{o}dinger equations are subcritical. For the critical case, one can consult  \cite{MR3018166,MR4182075,MR4393373} and references therein.

On the contrary, in our problem (\ref{e9}), no restriction
on the power $s$ is imposed: $s$ can be greater than $2\cdot2^{\ast}$, in this
case the nonlinear term $h(x)\left\vert u\right\vert ^{s-2}u$ is supercriticsal. There are also a few papers about supercritical problems, see e.g.\ \cite{MR3412303,MR3490340}. To study supercritical problems, one applies variational methods to get solutions of the subcritcal problem obtained by modifying $g(x,u)$ for $\vert u\vert $ large, then perform $L^\infty$-estimate to show that the solutions for the truncated problem have small $L^\infty$-norm, therefore  they are solutions of the original problem. We can see that $L^\infty$-estimate is a crucial step for this approach. Our approach for getting solutions of (\ref{e9}) does not require truncation and $L^\infty$-estimate.

Removing the quasilinear term $u\Delta(u^{2})$ from (\ref{e9}), our equation
reduces to a semilinear elliptic equation on $\mathbb{R}^{N}$, which was
studied by Tonkes \cite{MR1742223}. However, in \cite{MR1742223} $s\leq
2^{\ast}$ is still required. The result of \cite{MR1742223} was extended by
Liu--Li \cite{MR1971024}, where a corresponding $p$-Laplacian problem is
considered, and $s$ can be greater than the critical Sobolev exponent
$p^{\ast}$. Naturally, our strategy to study the supercritical quasilinear
problem (\ref{e9}) is motivated by \cite{MR1971024}. However, due to the above
mentioned nonlinear transformation $f$, more delicate analysis is needed. We
will see that the geometric properties of $f$ play an essential role in our investigation.

To state our result, for $p\in\left(  1,\infty\right)  $ we denote%
\[
p_{0}=\frac{2N}{2N-p\left(  N-2\right)  }\text{,\qquad}p^{\prime}=\frac
{p}{p-1}\text{.}%
\]
Note that $p^{\prime}$ is the H\"{o}lder conjugate exponent of $p$.

\begin{thm}
\label{t1}Assume that

\begin{enumerate}
\item[$\left(  k\right)  $] $k\in L^{q_{0}}(\mathbb{R}^{N})\cap L^{\infty}(\mathbb{R}^{N})$, $k\geq0$,
$k\not \equiv 0$,

\item[$\left(  h\right)  $] $h\in L^{1}(\mathbb{R}^{N})\cap L^{\infty}(\mathbb{R}^{N})$, $h\geq0$
\end{enumerate}
are satisfied, then the problem \eqref{e9} has a sequence of solutions
$\left\{  u_{n}\right\}  $ such that the energy $J(u_{n})<0$ and
$J(u_{n})\rightarrow0$ as $n\rightarrow\infty$.
\end{thm}

\begin{thm}
\label{t2}Under the assumptions $\left(  k\right)  $ and $\left(  h\right)  $,
the problem \eqref{e9} has a \emph{nonegative} solution $u$ such that $J(u)<0$.
\end{thm}

\begin{rem}
Our Theorem \ref{t2} is closely related to \cite{MR3250499}, where for $4\leq q<s<\infty$, the following problem%
\[
-\Delta u-\Delta(u^{2})=\lambda u+k(x)\left\vert u\right\vert ^{q-2}%
u-h(x)\left\vert u\right\vert ^{s-2}u\text{,\qquad}u\in H_{0}^{1}(\Omega)
\]
on a bounded domain $\Omega$ is considered; for $\lambda\in(\lambda^*,\bar{\lambda})$, a
nonnegative solution (at negative energy level) is obtained by the Ekeland
variational principle and the sub-super solution method.
\end{rem}

The paper is organized as follows. In Section 2 we review the definition of
the transformation $f$ and present some of its properties which are needed in
this paper. Since the exponent $s$ in (\ref{e9}) can be greater than the
critical Sobolev exponent $2\cdot2^{\ast}$, instead of the usual Sobolev
spaces $H^{1}(\mathbb{R}^{N})$ or $D^{1,2}(\mathbb{R}^{N})$, we introduce a
new space $E$ as the foundation of our functional framework. In Section 3 we investigate the
geometry and compactness of our energy functional $\Phi:E\rightarrow
\mathbb{R}$ and prove our theorems
via minimization argument and a variant of Clark's theorem proved in \cite{MR1828946}.

\section{Variational framework}

Following \cite{MR2029068,MR1949452}, we make the change of variables by
$u=f(v)$, where $f$ is an odd function defined by%
\begin{equation}
f^{\prime}(t)=\frac{1}{\sqrt{1+2f^{2}(t)}}\text{,\qquad}f(0)=0 \label{ode}%
\end{equation}
on $[0,+\infty)$. The proof of the following proposition can be found in
\cite{MR3121525,MR2029068} (some of them are obvious from (\ref{ode})).

\begin{prop}
\label{p1}The function $f$ possesses the following properties:

\begin{enumerate}
\item $f\in C^{\infty}(\mathbb{R})$ is strictly increasing, therefore is invertible.

\item $\left\vert f(t)\right\vert \leq\left\vert t\right\vert $, $f^{\prime
}(0)=1$, $\left\vert f^{\prime}(t)\right\vert \leq1$ for all $t\in\mathbb{R}$.

\item $\left\vert f(t)f^{\prime}(t)\right\vert \leq1$, $\left\vert
f(t)\right\vert \leq2^{1/4}\left\vert t\right\vert ^{1/2}$.

\item There exists a positive constant $\mu$ such that%
\begin{equation}
\left\vert f(t)\right\vert \geq\mu\left\vert t\right\vert \text{\quad for
}\left\vert t\right\vert \leq1\text{,\qquad}\left\vert f(t)\right\vert \geq
\mu\left\vert t\right\vert ^{1/2}\text{\quad for }\left\vert t\right\vert
\geq1\text{.} \label{mu}%
\end{equation}

\item For all $t\in\mathbb{R}$ we have $f^{2}(t)\geq f(t)f^{\prime}(t)t\ge\dfrac12f^2(t)$.
\end{enumerate}
\end{prop}

Motivated by \cite{MR1971024}, let $E$ be the completion of $C_{0}^{\infty}(\mathbb{R}^{N})$
under the norm
\begin{align}
\left\Vert v\right\Vert  &  =\left\Vert v\right\Vert _{D}+\vert h^{2/s}%
v\vert_{s/2}\nonumber\\
&  =\left(  \int\left\vert \nabla v\right\vert ^{2}\right)  ^{1/2}+\left(
\int h\left\vert v\right\vert ^{s/2}\right)  ^{2/s}\text{,} \label{n}%
\end{align}
where $\left\Vert \cdot\right\Vert _{D}$ and $\left\vert \,\cdot\,\right\vert
_{p}$ are the standard $D^{1,2}$-norm and $L^{p}$-norm ($p\in\left[
1,\infty\right]  $), respectively. Note that if following \cite{MR1971024}
directly, one may tend to define the norm as%
\[
\left\Vert v\right\Vert =\left\Vert v\right\Vert _{D}+\vert h^{1/s}v\vert _{s}\text{.}%
\]
Our definition (\ref{n}) takes the structure of (\ref{e9}) and the growth property $\left\vert f(t)\right\vert
\leq c\left\vert t\right\vert ^{1/2}$ of $f$ into account. It turns out
that this is the correct choice.

\begin{rem}
When $h\equiv0$, our space $E$ reduces to the standard Sobolev space $D^{1,2}(\mathbb{R}^N)$.
\end{rem}

To present the variational framework for our argument, we need the following lemma.

\begin{lem}
\label{l0}If $\phi\in C_{0}^{\infty}(\mathbb{R}^{N})$, then
\begin{equation}
\xi=\frac{\phi}{f^{\prime}(v)}=\sqrt{1+2f^{2}(v)}\phi\label{el}%
\end{equation}
belongs to $E$.
\end{lem}

\begin{pf}
Take $R>0$ such that $\operatorname*{supp}\phi\subset B_{R}$, where $B_R$ is the $R$-ball in $\mathbb{R}^N$. Since
\[
\left(  1+2f^{2}(v)\right)  ^{s/4}\leq C\big(  1+\left\vert v\right\vert
^{s/2}\big)  \text{,}%
\]
we have%
\begin{align}
\int h\left\vert \xi\right\vert ^{s/2}  &  =\int h\big(  1+2f^{2}(v)\big)
^{s/4}\left\vert \phi\right\vert ^{s/2}\nonumber\\
&  \leq C\left\vert \phi\right\vert _{\infty}^{s/2}\int h\big(  1+\left\vert
v\right\vert ^{s/2}\big) \nonumber\\
&  \leq C\left\vert \phi\right\vert _{\infty}^{s/2}\left(  \left\vert
h\right\vert _{1}+\int h\left\vert v\right\vert ^{s/2}\right)  <\infty\text{.}
\label{E1}%
\end{align}
Now we estimate the $D^{1,2}$-norm of $\xi$. Because $v\in D^{1,2}(\mathbb{R}^{N})$, we have
$v\in L_{\mathrm{loc}}^{2}(\mathbb{R}^{N})$, therefore%
\begin{align}
\int\left\vert \nabla\xi\right\vert ^{2}  &  =\int_{B_{R}}\left\vert \left(
1+2f^{2}(v)\right)  ^{1/2}\nabla\phi+\frac{2f(v)f^{\prime}(v)\phi}%
{\sqrt{1+2f^{2}(v)}}\nabla v\right\vert ^{2}\nonumber\\
&  \leq\int_{B_{R}}\left(  \left(  1+2v^{2}\right)  \left\vert \nabla
\phi\right\vert ^{2}+4\left\vert v\right\vert \left\vert \phi\right\vert
\left\vert \nabla\phi\right\vert \left\vert \nabla v\right\vert +4\phi
^{2}\left\vert \nabla v\right\vert ^{2}\right) \nonumber\\
&  \leq m\int_{B_{R}}\left(  \left(  1+2v^{2}\right)  +4\left\vert
v\right\vert \left\vert \nabla v\right\vert +4\left\vert \nabla v\right\vert
^{2}\right)  <\infty\text{,} \label{E2}%
\end{align}
where $m=\left(\left\vert \phi\right\vert _{\infty}+\left\vert \nabla
\phi\right\vert _{\infty}\right)  ^{2}$. Combining (\ref{E1}) and (\ref{E2})
we see that $\xi\in E$.
\end{pf}

By the growth properties of $f$, it is easy to see that under our assumptions
on $k$ and $h$, the functional%
\[
\Phi(v)=J(f(v))=\frac{1}{2}\int\left\vert \nabla v\right\vert ^{2}-\frac{1}%
{q}\int k\left\vert f(v)\right\vert ^{q}+\frac{1}{s}\int h\left\vert
f(v)\right\vert ^{s}
\]
is well defined and of class $C^{1}$ on the Banach space $E$, with derivative
given by%
\begin{align*}
\langle\Phi^{\prime}(v),\xi\rangle=\int\nabla v\cdot\nabla\xi &  -\int
k\left\vert f(v)\right\vert ^{q-2}f(v)f^{\prime}(v)\xi\\
& \qquad +\int h\left\vert f(v)\right\vert ^{s-2}f(v)f^{\prime}(v)\xi
\end{align*}
for $v,\xi\in E$. Moreover, if $v$ is a critical point of $\Phi:E\rightarrow
\mathbb{R}$, by Lemma \ref{l0} for $\phi\in C_{0}^{\infty}(\mathbb{R}^{N})$ we
have $\xi=\phi/f^{\prime}(v)\in E$. Hence $\left\langle \Phi^{\prime}%
(v),\xi\right\rangle =0$ and from which let $u=f(v)$, by standard computation
we get%
\[
0=\left.  \frac{\mathrm{d}}{\mathrm{d}t}\right\vert _{t=0}J(u+t\phi)\text{,}%
\]
which means that $u$ is a weak solution of the problem (\ref{e9}). Therefore,
to prove our theorems it suffices to find critical points of $\Phi:E\to\mathbb{R}$.
This is the task of the next section.

The following proposition justifies our effort to find solutions of (\ref{e9})
at negative energy levels.

\begin{prop}
\label{pp}Suppose $s\ge4$. If $v\in E$ is a critical point of $\Phi$, then
$\Phi(v)\leq0$.
\end{prop}

\begin{pf}
Let $c=\Phi(v)$, that is%
\begin{equation}
2c=\int\left\vert \nabla v\right\vert ^{2}-\frac{2}{q}\int k\left\vert
f(v)\right\vert ^{q}+\frac{2}{s}\int h\left\vert f(v)\right\vert ^{s}\text{.}
\label{ew}%
\end{equation}
Testing $\Phi^{\prime}(v)$ by $v$, we have%
\begin{align}
0  &  =\langle\Phi^{\prime}(v),v\rangle\nonumber\\
&  =\int\left\vert \nabla v\right\vert ^{2}-\int k\left\vert f(v)\right\vert
^{q-2}f(v)f^{\prime}(v)v+\int h\left\vert f(v)\right\vert ^{s-2}f(v)f^{\prime
}(v)v\text{.} \label{er}%
\end{align}
Since $q<2$, $s\ge4$, from 
(\ref{ew}), (\ref{er}) and Proposition \ref{p1} (5), we obtain%
\begin{align*}
2c  &  =\int k\left\vert f(v)\right\vert ^{q-2}\left[  f(v)f^{\prime
}(v)v-\frac{2}{q}f^{2}(v)\right] \\
&  \qquad+\int h\left\vert f(v)\right\vert ^{s-2}\left[  \frac{2}{s}%
f^{2}(v)-f(v)f^{\prime}(v)v\right]  \leq0\text{,}%
\end{align*}
as desired.
\end{pf}

\begin{rem}
The idea of \cite{MR1971024}, which in turn was inspired by an unpublished
preprint \cite{FW}, has also been employed in \cite{MR3286966}, where a
supercritical Schr\"{o}dinger-Poisson system%
\[
\left\{
\begin{array}
[c]{ll}%
-\Delta u+a(x)u+\phi u=k(x)\left\vert u\right\vert ^{q-2}u-h(x)\left\vert
u\right\vert ^{s-2}u & x\in\mathbb{R}^{3}\text{,}\\
-\Delta\phi=u^{2}\text{, }\lim\limits_{\left\vert x\right\vert \rightarrow
\infty}\phi(x)=0 & x\in\mathbb{R}^{3}%
\end{array}
\right.
\]
is studied. In \cite{FW,MR1971024,MR3286966} no nonlinear transformation like
$f$ is involved, therefore our work is not a trivial application of the idea
in these papers.
\end{rem}

\section{Proof of Theorems \ref{t1} and \ref{t2}}

In this section, we will show that, under our assumptions, $\Phi
:E\rightarrow\mathbb{R}$ is coercive and satisfies the Palais-Smale condition,
then prove our theorems by minimization method and a variant of the classical Clark's theorem.

\begin{lem}\label{ll}
$\Phi:E\rightarrow\mathbb{R}$ is coercive.
\end{lem}

\begin{pf}
Let $\ell$ be the norm of the embedding $D^{1,2}(\mathbb{R}^{N})\hookrightarrow L^{2^{\ast}}(\mathbb{R}^{N})$.
Since $qq_{0}^{\prime}=2^{\ast}$, for $v\in E$ we have%
\begin{align}
\int k\left\vert v\right\vert ^{q}  &  \leq\left\vert k\right\vert _{q_{0}%
}\left\vert \left\vert v\right\vert ^{q}\right\vert _{q_{0}^{\prime}%
}\nonumber\\
&  =\left\vert k\right\vert _{q_{0}}\left\vert v\right\vert _{2^{\ast}}%
^{q}\leq\ell^{q}\left\vert k\right\vert _{q_{0}}\left\Vert v\right\Vert
_{D}^{q}\text{.} \label{qq}%
\end{align}
If $\Phi$ is not coercive, there is a sequence $\left\{  v_{n}\right\}  $ in
$E$ such that $\sup_{n}\Phi(v_{n})<+\infty$ and%
\begin{align}
\left\Vert v_{n}\right\Vert =\left\Vert v_{n}\right\Vert _{D}+\vert h^{2/s}%
v_{n}\vert_{s/2}\rightarrow+\infty\text{.} \label{xy}%
\end{align}
We claim that $\{v_n\}$ is bounded in $D^{1,2}(\mathbb{R}^{N})$.  Otherwise along a subsequence  we have $\Vert
v_{n_{i}}\Vert _{D}\rightarrow\infty$, using (\ref{qq}) and noting $q<2$ we have
\begin{align*}
\Phi(v_{n_{i}})  &  \geq\frac{1}{2}\int\vert \nabla v_{n_{i}}\vert
^{2}-\frac{1}{q}\int k\vert v_{n_{i}}\vert ^{q}+\frac{1}{s}\int
h\vert f(v_{n_{i}})\vert ^{s}\\
&  \geq\frac{1}{2}\Vert v_{n_{i}}\Vert _{D}^{2}-\frac{1}{q}\ell
^{q}\left\vert k\right\vert _{q_{0}}\Vert v_{n_{i}}\Vert _{D}%
^{q}+\frac{1}{s}\int h\vert f(v_{n_{i}})\vert ^{s}\\
&  \geq\frac{1}{2}\Vert v_{n_{i}}\Vert _{D}^{2}-\frac{1}{q}\ell
^{q}\left\vert k\right\vert _{q_{0}}\Vert v_{n_{i}}\Vert _{D}%
^{q}\rightarrow+\infty\text{,}%
\end{align*}
contradicting $\sup_{n}\Phi(v_{n})<+\infty$. Therefore $\sup_{n}\left\Vert
v_{n}\right\Vert _{D}<\infty$ and from (\ref{xy}) we have
\[
\int h\left\vert v_{n}\right\vert ^{2/s}=\vert h^{2/s}v_{n}\vert _{s/2}^{s/2}%
\rightarrow+\infty\text{.}%
\]
Using (\ref{mu}) we get%
\begin{align*}
\int h\left\vert f(v_{n})\right\vert ^{s}  &  =\int_{\left\vert v_{n}%
\right\vert \leq1}h\left\vert f(v_{n})\right\vert ^{s}+\int_{v_{n}%
>1}h\left\vert f(v_{n})\right\vert ^{s}\\
&  \geq\mu\int_{v_{n}>1}h\left\vert v_{n}\right\vert ^{s/2}=\mu\int
h\left\vert v_{n}\right\vert ^{2/s}-\mu\int_{\left\vert v_{n}\right\vert
\leq1}h\left\vert v_{n}\right\vert ^{2/s}\\
&  \geq\mu\int h\left\vert v_{n}\right\vert ^{2/s}-\mu\left\vert h\right\vert
_{1}\rightarrow+\infty\text{.}%
\end{align*}
Since $q<2$, we end up at a contradiction:%
\begin{align*}
\Phi(v_{n})    \geq\frac{1}{2}\left\Vert v_{n}\right\Vert _{D}^{2}&-\frac
{1}{q}\ell^{q}\left\vert k\right\vert _{q_{0}}\left\Vert v_{n}\right\Vert_D
^{q}\\
&+\frac{1}{s}\int h\left\vert f(v_{n})\right\vert ^{s}
\rightarrow+\infty\text{.}%
\end{align*}
The proof is completed.
\end{pf}

\begin{lem}
\label{l1}Given $a\in\mathbb{R}$, the function $\eta:\mathbb{R}\rightarrow
\mathbb{R}$, $\eta(t)=\left\vert f(t)\right\vert ^{s}$, is convex. Hence for
$\alpha,\beta\in\mathbb{R}$ we have%
\begin{equation}
\left\vert f(\alpha)\right\vert ^{s}\leq\left\vert f(\beta)\right\vert
^{s}+s\left\vert f(\alpha)\right\vert ^{s-2}f(\alpha)f^{\prime}(\alpha)\left(
\alpha-\beta\right)  \text{.} \label{ea}%
\end{equation}

\end{lem}

\begin{pf}
Obviously $\eta$ is smooth and even. For $t\geq0$, because $s>2$, using
(\ref{ode}) we have
\begin{align*}
\eta^{\prime}  &  =sf^{s-1}f^{\prime}=\frac{sf^{s-1}}{\sqrt{1+f^{2}}}%
\text{,}\\
\frac{\eta^{\prime\prime} }s  &  =\left(  \frac{f^{s-1}}{\sqrt{1+f^{2}}%
}\right)  ^{\prime}\\
&  =\frac{\left(  s-1\right)  f^{s-2}f^{\prime}\sqrt{1+f^{2}}-f^{s-1}%
\dfrac{ff^{\prime}}{\sqrt{1+f^{2}}}}{1+f^{2}}\\
&  =\frac{\left(  s-2\right)  f^{s}f^{\prime}+\left(  s-1\right)
f^{s-2}f^{\prime}}{\left(  1+f^{2}\right)  ^{3/2}}\geq0\text{.}%
\end{align*}
Because $\eta^{\prime\prime}$ is also even, we see that $\eta^{\prime\prime
}(t)\geq0$ for all $t\in\mathbb{R}$, and $\eta$ is convex.
\end{pf}

\begin{lem}
$\Phi$ satisfies the Palais-Smale condition.
\end{lem}

\begin{pf}
Let $\left\{  v_{n}\right\}  \subset E$ be a $\left(  PS\right)  $ sequence.
Lemma \ref{l1} implies that $\left\{  v_{n}\right\}  $ is bounded. Thus $\left\{  v_{n}\right\}  $
and $\big\{  h^{2/s}v_{n}\big\}  $ are bounded in $D^{1,2}(\mathbb{R}^{N})$ and $L^{s/2}(\mathbb{R}^{N})$,
respectively. Up to a subsequence we have%
\begin{equation}
v_{n}\rightharpoonup v\text{\quad in }D^{1,2}(\mathbb{R}^{N})\text{,}\qquad h^{2/s}%
v_{n}\rightharpoonup h^{2/s}v\text{\quad in }L^{s/2}(\mathbb{R}^{N})\text{.} \label{ee}%
\end{equation}
From this it is clear that $v\in E$. Moreover, according to \cite[Lemma 1]{MR1457116}, our condition $\left(
k\right)  $ implies that the functional
\begin{align}\psi:D^{1,2}(\mathbb{R}^{N})\to\mathbb{R}\text{,}\qquad \psi(v)=\int k\left\vert
v\right\vert ^{q}\label{kq}
\end{align}
is weakly continuous on $D^{1,2}(\mathbb{R}^{N})$, thus%
\begin{align}
\int k\left\vert v_{n}-v\right\vert ^{q}\rightarrow0\text{.} \label{k}%
\end{align}
Firstly we want to show%
\begin{equation}
\langle\Phi^{\prime}(v),v_{n}-v\rangle\rightarrow0\text{.} \label{e0}%
\end{equation}
Since we don't know whether our space $E$ is reflexive, we could not get
$v_{n}\rightharpoonup v$ in $E$ and deduce (\ref{e0}). Therefore we adapt the
following argument.

Because $\left\vert f(t)\right\vert \leq\left\vert t\right\vert $ and
$\left\vert f^{\prime}(t)\right\vert \leq1$ (see Proposition \ref{p1}), using
H\"{o}lder inequality and (\ref{k}) we get%
\begin{align}
&  \left\vert \int k\left\vert f(v)\right\vert ^{q-2}f(v)f^{\prime}(v)\left(
v_{n}-v\right)  \right\vert \leq\int k\left\vert v\right\vert ^{q-1}\left\vert
v_{n}-v\right\vert \nonumber\\
&  \qquad\leq\left(  \int k\left\vert v_{n}-v\right\vert ^{q}\right)
^{1/q}\left(  \int k\left\vert v\right\vert ^{q}\right)  ^{\left(  q-1\right)
/q}\rightarrow0\text{.} \label{x}%
\end{align}
Since $\left\vert f(t)f^{\prime}(t)\right\vert \leq1$ and $\left\vert
f(t)\right\vert \leq2^{1/4}\left\vert t\right\vert ^{1/2}$, noting $\left(
s/2\right)  ^{\prime}=s/\left(  s-2\right)  $ we have%
\begin{align*}
\int\left\vert h^{1-2/s}\left\vert f(v)\right\vert ^{s-2}f(v)f^{\prime
}(v)\right\vert ^{\left(  s/2\right)  ^{\prime}}  &  \leq\int h\left(
\left\vert f(v)\right\vert ^{s-2}\right)  ^{s/\left(  s-2\right)  }\\
&  =\int h\left\vert f(v)\right\vert ^{s}\leq2^{s/4}\int h\left\vert
v\right\vert ^{s/2}<\infty\text{,}%
\end{align*}
that is $h^{1-2/s}\left\vert f(v)\right\vert ^{s-2}f(v)f^{\prime}(v)\in
L^{\left(  s/2\right)  ^{\prime}}(\mathbb{R}^{N})$. Using $h^{2/s}v_{n}\rightharpoonup
h^{2/s}v$ in $L^{s/2}(\mathbb{R}^{N})$ we get%
\begin{align}
&  \int h\left\vert f(v)\right\vert ^{s-2}f(v)f^{\prime}(v)\left(
v_{n}-v\right) \nonumber\\
&  \qquad=\int h^{2/s}\left(  v_{n}-v\right)  \cdot h^{1-2/s}\left\vert
f(v)\right\vert ^{s-2}f(v)f^{\prime}(v)\rightarrow0\text{.} \label{e2}%
\end{align}
Combining (\ref{e2}) with (\ref{ee}) and (\ref{x}) we get%
\begin{align*}
\langle\Phi^{\prime}(v),v_{n}-v\rangle &  =\int\nabla v\cdot\nabla\left(
v_{n}-v\right)  -\int k\left\vert f(v)\right\vert ^{q-2}f(v)f^{\prime
}(v)\left(  v_{n}-v\right) \\
&  \qquad\qquad+\int h\left\vert f(v)\right\vert ^{s-2}f(v)f^{\prime
}(v)\left(  v_{n}-v\right)  \rightarrow0\text{,}%
\end{align*}
our claim (\ref{e0}) follows.

Next, using $\left\vert f(t)\right\vert \leq\left\vert t\right\vert $ and
$\left\vert f^{\prime}(t)\right\vert \leq1$ again, for
\[
\Omega_{n}:=\left\vert f(v_{n})\right\vert ^{q-2}f(v_{n})f^{\prime}%
(v_{n})-\left\vert f(v)\right\vert ^{q-2}f(v)f^{\prime}(v)\text{,}
\]
we have%
\begin{align*}
\left\vert \Omega_{n}\right\vert    \leq\left\vert f(v_{n})\right\vert
^{q-1}+\left\vert f(v)\right\vert ^{q-1}\leq\left\vert v_{n}\right\vert
^{q-1}+\left\vert v\right\vert ^{q-1}\text{.}
\end{align*}
Since $\{k^{1/q}v_{n}\}$ is bounded in $L^{q}(\mathbb{R}^{N})$, 
\begin{align*}
\left\vert \int k\left\vert \Omega_{n}\right\vert ^{q^{\prime}}\right\vert  &
\leq2^{q^{\prime}}\left(  \int k\left(  \left\vert v_{n}\right\vert
^{q}+\left\vert v\right\vert ^{q}\right)  \right) \\
&  \leq2^{q^{\prime}+1}\sup_{n}\vert k^{1/q}v_{n}\vert_{q}^{q}=:M<\infty\text{.}%
\end{align*}
Thus using H\"{o}lder inequality and (\ref{k}) we deduce%
\begin{align}
&  \left\vert \int k\left(  \left\vert f(v_{n})\right\vert ^{q-2}%
f(v_{n})f^{\prime}(v_{n})-\left\vert f(v)\right\vert ^{q-2}f(v)f^{\prime
}(v)\right)  \left(  v_{n}-v\right)  \right\vert \nonumber\\
&  \qquad\leq\int k^{1/q}\left\vert v_{n}-v\right\vert \cdot k^{1/q^{\prime}%
}\vert\Omega_{n}\vert\nonumber\\
&  \qquad\leq\left(  \int k\left\vert v_{n}-v\right\vert ^{q}\right)
^{1/q}\left(  \int k\left\vert \Omega_{n}\right\vert ^{q^{\prime}}\right)
^{1/q^{\prime}}\nonumber\\
&  \qquad\leq M^{1/q^{\prime}}\left(  \int k\left\vert v_{n}-v\right\vert
^{q}\right)  ^{1/q}\rightarrow0\text{.} \label{e1}%
\end{align}
On the other hand, noting that the function%
\[
t\mapsto s\left\vert f(t)\right\vert ^{s-2}f(t)f^{\prime}(t)
\]
is increasing (it is the derivative of the convex function $\eta$ given in
Lemma \ref{l1}), we get%
\[
H_{n}:=\int h\left(  \left\vert f(v_{n})\right\vert ^{s-2}f(v_{n})f^{\prime
}(v_{n})-\left\vert f(v)\right\vert ^{s-2}f(v)f^{\prime}(v)\right)  \left(
v_{n}-v\right)  \geq0\text{.}%
\]
Now, using (\ref{e0}) and (\ref{e1}) we get%
\begin{align}
o(1)  &  =\langle\Phi^{\prime}(v_{n})-\Phi^{\prime}(v),v_{n}-v\rangle
\nonumber\\
&  =\int\left\vert \nabla(v_{n}-v)\right\vert ^{2}\nonumber\\
&  \qquad-\int k\left(  \left\vert f(v_{n})\right\vert ^{q-2}f(v_{n}%
)f^{\prime}(v_{n})-\left\vert f(v)\right\vert ^{q-2}f(v)f^{\prime}(v)\right)
\left(  v_{n}-v\right) \nonumber\\
&  \qquad+\int h\left(  \left\vert f(v_{n})\right\vert ^{s-2}f(v_{n}%
)f^{\prime}(v_{n})-\left\vert f(v)\right\vert ^{s-2}f(v)f^{\prime}(v)\right)
\left(  v_{n}-v\right) \nonumber\\
&  =\int\left\vert \nabla\left(  v_{n}-v\right)  \right\vert ^{2}%
+H_{n}+o(1)\text{.} \label{g}%
\end{align}
Consequently, noting $H_n\ge0$ we deduce
\begin{equation}
v_{n}\rightarrow v\text{\quad in }D^{1,2}(\mathbb{R}^{N})\text{,\qquad}H_{n}\rightarrow
0\text{.} \label{e}%
\end{equation}
Since $H_{n}\rightarrow0$, from (\ref{e2}) we have%
\[
\int h\left\vert f(v_{n})\right\vert ^{s-2}f(v_{n})f^{\prime}(v_{n})\left(
v_{n}-v\right)  \rightarrow0\text{.}%
\]
Replacing $\alpha$ and $\beta$ in (\ref{ea}) with $v_{n}$ and $v$
respectively, we get%
\begin{align*}
\varlimsup_{n\rightarrow\infty}\int h\left\vert f(v_{n})\right\vert ^{s}  &
\leq\int h\left\vert f(v)\right\vert ^{s}+s\lim_{n\rightarrow\infty}\int
h\left\vert f(v_{n})\right\vert ^{s-2}f(v_{n})f^{\prime}(v_{n})\left(
v_{n}-v\right) \\
&  =\int h\left\vert f(v)\right\vert ^{s}\text{.}%
\end{align*}
Combining this with the easy consequence
\[
\int h\left\vert f(v)\right\vert ^{s}\leq\varliminf_{n\rightarrow\infty}\int
h\left\vert f(v_{n})\right\vert ^{s}%
\]
of $v_{n}\rightarrow v$ a.e. in $\mathbb{R}^{N}$ and Fatou'e lemma, we get%
\begin{equation}
\int h\left\vert f(v_{n})\right\vert ^{s}\rightarrow\int h\left\vert
f(v)\right\vert ^{s}\text{.} \label{eh}%
\end{equation}
Now, noting the following consequence of (\ref{mu}):
\[
h\left\vert v_{n}\right\vert ^{s/2}\leq h+\frac{1}{\mu^s}h\left\vert
f(v_{n})\right\vert ^{s}%
\]
and $h\left\vert v_{n}\right\vert ^{s/2}\rightarrow h\left\vert v\right\vert
^{s/2}$ a.e.\ in $\mathbb{R}^{N}$, by the generalized Lebesgue dominating
theorem (see Proposition \ref{pl} below) and (\ref{eh}) we get%
\[
\int h\left\vert v_{n}\right\vert ^{s/2}%
\rightarrow\int h\left\vert v\right\vert ^{s/2}\text{.}%
\]
That is to say $\vert h^{2/s}v_{n}\vert _{s/2}\to \vert h^{2/s}v\vert _{s/2}$. But $h^{2/s}v_{n}\rightharpoonup h^{2/s}v$ in $L^{s/2}(\mathbb{R}^{N})$, we deduce
$h^{2/s}v_{n}\rightarrow h^{2/s}v$ in $L^{s/2}(\mathbb{R}^{N})$. Combining this with (\ref{e})
we get%
\begin{align*}
\left\Vert v_{n}-v\right\Vert  &  =\left(  \int\left\vert \nabla\left(
v_{n}-v\right)  \right\vert ^{2}\right)  ^{1/2}+\left(  \int h\left\vert
v_{n}-v\right\vert ^{s/2}\right)  ^{2/s}\\
&  =\left\Vert v_{n}-v\right\Vert _{D}+\big\vert h^{2/s}v_{n}-h^{2/s}%
v\big\vert _{s/2}\rightarrow0\text{.}%
\end{align*}
Thus $v_{n}\rightarrow v$ in $E$.
\end{pf}

For the reader's convenience, we quote the generalized Lebesgue dominating
theorem as follow.

\begin{prop}
\label{pl} Let $f_{n},g_{n}:\Omega\rightarrow\mathbb{R}$ be measurable
functions over the measurable set $\Omega$, $f_{n}\rightarrow f$ a.e.\ in
$\Omega$, $g_{n}\rightarrow g$ a.e.\ in $\Omega$, $\left\vert f_{n}\right\vert
\leq g_{n}$. Then%
\[
\int_{\Omega}\left\vert f_{n}-f\right\vert \rightarrow0
\]
provided $\int_{\Omega}g_{n}\rightarrow\int_{\Omega}g$ and $\int_{\Omega
}g<+\infty$.
\end{prop}

\begin{rem}
As is well known, to prove Proposition \ref{pl} we apply Fatou's lemma to
\[
F_n:=g_n+g-\vert f_n-f\vert \text{.}
\]
When $g_n=g$ does not depend on $n$, Proposition \ref{pl} reduces to the usual Lebesgue dominating theorem. 
\end{rem}

Having verified the $\left(  PS\right)  $ condition, we need the following
variant of Clark's theorem (see \cite{MR296777} or \cite[Theorem
9.1]{MR845785} for the classical Clark's theorem) to produce the desired solutions of our problem (\ref{e9}).

\begin{prop}
[{\cite[Lemma 2.4]{MR1828946}}]\label{p}Let $E$ be a Banach space and $\Phi\in
C^{1}(E,\mathbb{R})$ be an even coercive functional satisfying the $\left(
PS\right)  $ condition and $\Phi(0)=0$. If for any $n\in\mathbb{N}$, there is
an $n$-dimensional subspace $X_{n}$ and $\rho_{n}>0$ such that%
\[
\sup_{X_{n}\cap S_{\rho_{_n}}}\Phi<0\text{,}%
\]
where $S_{r}=\left\{  u\in E\vert \,\left\Vert u\right\Vert =r\right\}$, then
$\Phi$ has a sequence of critical values $c_{n}<0$ satisfying $c_{n}%
\rightarrow0$.
\end{prop}

\begin{pf}
[Proof of Theorem \ref{t1}]
Given $n\in\mathbb{N}$, let $X_{n}$ be an $n$-dimensional subspace of $X$. Since the norms $\left\Vert \cdot\right\Vert $ and $\left\vert \,\cdot\,
\right\vert _{\infty}$ are equivalent on $X_{n}$, there is $\vartheta>0$ such
that%
\[
\left\vert v\right\vert _{\infty}\leq\vartheta\left\Vert v\right\Vert
\qquad\text{for all }v\in X_{n}\text{.}%
\]
Because $h\in L^{1}(\mathbb{R}^{N})$, we have
\[
\left\vert \int h\left\vert v\right\vert ^{s}\right\vert \leq\left\vert
v\right\vert _{\infty}^{s}\int h\leq\vartheta^{s}\left\Vert v\right\Vert ^{s}\vert
h\vert_1<\infty\text{.}%
\]
Thus we have a well-defined $s$-homogeneous functional $H:X_{n}\rightarrow\mathbb{R}%
$,%
\[
H(v)=\int h\left\vert v\right\vert ^{s}\text{.}%
\]
Using the Lebesgue dominating
theorem, it is easy to see that $H$ is continuous.

Since $f^{\prime}(0)=1$, there is $\delta\in\left(  0,1\right)  $ such that%
\begin{equation}
\frac{1}{2}\left\vert t\right\vert \leq\left\vert f(t)\right\vert
\leq\left\vert t\right\vert \text{,\qquad for }t\in\left[  -\delta
,\delta\right]  \text{.} \label{f}%
\end{equation}
Because $\dim X_{n}<\infty$, the compactness of
\[
X_{n}\cap S_{1}=\left\{  v\in X_{n}\vert\,\left\Vert v\right\Vert =1\right\}
\]
and $k\not \equiv 0$ implies that both%
\begin{equation}
A=\inf_{\varphi\in X_{n}\cap S_{1}}\int k\left\vert \varphi\right\vert
^{q}\text{\qquad and\qquad}B=\sup_{\varphi\in X_{n}\cap S_{1}}H(\varphi)=\sup_{\varphi\in X_{n}\cap S_{1}}\int h\vert \varphi\vert^s\label{s}
\end{equation}
are finite, $A>0$, $B\ge0$. Since the norms $\left\Vert \cdot\right\Vert $
and $\left\vert \,\cdot\,
\right\vert _{\infty}$ on $X_{n}$ are equivalent, noting $q<2<s$, we can choose
$\rho_{n}>0$ such that if $v\in X_{n}$, $\left\Vert v\right\Vert =\rho_{n}$,
then $\left\vert v\right\vert _{\infty}\leq\delta$ and%
\begin{align}
\theta_{n}:=\frac{1}{2}\rho_{n}^{2}-\frac{A}{2^{q}q}\rho_{n}^{q}+\frac{
B}{s}\rho_{n}^{s}<0\text{.}%
\label{qx}
\end{align}
Now using (\ref{f}) and (\ref{s}), for $v\in X_{n}\cap
S_{\rho_{_n}}$ we have $\left\vert v\right\vert _{\infty}\leq\delta$ and
\begin{align}
\Phi(v)  &  =\frac{1}{2}\int\left\vert \nabla v\right\vert ^{2}-\frac{1}%
{q}\int k\left\vert f(v)\right\vert ^{q}+\frac{1}{s}\int h\left\vert
f(v)\right\vert ^{s}\nonumber\\
&  \leq\frac{1}{2}\left\Vert v\right\Vert _{D}^{2}-\frac{1}{2^{q}q}\int
k\left\vert v\right\vert ^{q}+\frac{1}{s}\int h\left\vert v\right\vert
^{s}\nonumber\\
&  \leq\frac{1}{2}\left\Vert v\right\Vert ^{2}-\frac{A}{2^{q}q}\left\Vert
v\right\Vert ^{q}+\frac{B}{s}\left\Vert v\right\Vert ^{s}=\theta_n\text{.}
\label{xx}%
\end{align}
From this, using (\ref{qx}) it is clear that%
\begin{align}
\sup_{X_{n}\cap S_{\rho_{_n}}}\Phi\leq\theta_{n}<0\text{.}%
\label{n}
\end{align}
Since our $\Phi$ is an even coercive functional satisfying the $\left(
PS\right)  $ condition and $\Phi(0)=0$, applying Proposition \ref{p} we know
that $\Phi$ has a sequence of critical points $\left\{  v_{n}\right\}  $ such
that%
\[
J(u_{n})=\Phi(v_{n})<0\text{,\qquad}J(u_{n})\rightarrow0\text{,}%
\]
where $u_{n}=f(v_{n})$ are the desired solutions of (\ref{e9}).
\end{pf}

\begin{pf}
[Proof of Theorem \ref{t2}]We know that $\Phi$ is bounded from below. Since $\Phi
(v_{n})=\Phi(\left\vert v_{n}\right\vert )$, we may take a minimization
sequence $\left\{  v_{n}\right\}  $ such that $v_{n}\geq0$ and%
\[
\Phi(v_{n})\rightarrow c:=\inf_{E}\Phi\text{,}
\]
where $c<0$ because from (\ref{n}) we know that $\Phi$ can take negative values. By Lemma \ref{ll} we know that $\left\{  v_{n}\right\}  $ is bounded in $D^{1,2}(\mathbb{R}^{N})$. Thus we may assume
\[
v_{n}\rightharpoonup v\text{\quad in }D^{1,2}(\mathbb{R}^{N})\text{,}\qquad v_n\to v\text{\quad a.e.\ in }\mathbb{R}^{N}
\]
for some nonnegative $v\in E$.

Since $h\vert f(v_n)\vert ^s\to h\vert f(v)\vert ^s$ a.e.\ in $\mathbb{R}^{N}$, using Fatou's lemma we get
\begin{align}
\int\vert \nabla v\vert ^2\le\varliminf_{n\to\infty}\int\vert \nabla v_n\vert ^2\text{,}\qquad\int h\vert f(v)\vert ^s\le\varliminf_{n\to\infty}\int h\vert f(v_n)\vert ^s\text{.}\label{f1}
\end{align}
By the weak continuity of the functional $\psi$ defined in (\ref{kq}), from $v_{n}\rightharpoonup v$ in $D^{1,2}(\mathbb{R}^{N})$ we have
\[
\int k\vert v_n\vert^q\to\int k\vert v\vert^q\text{.}
\]
Since $k\vert f(v_n)\vert^q\le k\vert v_n\vert^q$, applying Proposition \ref{pl} we get
\begin{align}
\int k\vert f(v)\vert ^q=\lim_{n\to\infty}\int k\vert f(v_n)\vert ^q\text{.}\label{f2}
\end{align}
From (\ref{f1}) and (\ref{f2}) we get
\begin{align*}
\Phi(v)&=\frac{1}{2}\int\left\vert \nabla v\right\vert ^{2}+\frac{1}{s}\int h\left\vert
f(v)\right\vert ^{s}-\frac{1}%
{q}\int k\left\vert f(v)\right\vert ^{q}\\
&\le\varliminf_{n\to\infty}\left(\frac{1}{2}\int\left\vert \nabla v_n\right\vert ^{2}+\frac{1}{s}\int h\left\vert
f(v_n)\right\vert ^{s}-\frac{1}%
{q}\int k\left\vert f(v_n)\right\vert ^{q}\right)\\
&=\varliminf_{n\to\infty}\Phi(v_n)=c\text{.}
\end{align*}
Therefore $\Phi(v)=c$ and $v$ is a nonnegative critical point of $\Phi$.

Since $f(t)$ has the same sign as $t$, $u=f(v)$ is a nonnegative solution of (\ref{e9}) at negative energy level $J(u)=\Phi(v)=c$.
\end{pf}

\begin{rem}
Under the same assumptions on $k$ and $h$, similar results holds for%
\[
\left\{
\begin{array}
[c]{l}%
-\Delta u-u\Delta(u^{2})=\lambda g(x)u+k(x)\left\vert u\right\vert
^{q-2}u-h(x)\left\vert u\right\vert ^{s-2}u\text{,}\\
u\in D^{1,2}(\mathbb{R}^{N})
\end{array}
\right.
\]
if $\lambda\in(  \lambda_{1}^{-},\lambda_{1}^{+})  $, where $g\in
L^{N/2}(\mathbb{R}^{N})\cap L^{\infty}(\mathbb{R}^{N})$, $\lambda_{1}^{\pm}$ are the principle
positive/negative eigenvalues of $-\Delta u=\lambda g(x)u$ on $D^{1,2}(\mathbb{R}^{N})$; see
e.g.\ \cite{MR1390979} for discussion about this eigenvalue problem. The reason is that if $\lambda\in(  \lambda_{1}^{-},\lambda_{1}^{+})  $ then there is $\kappa>0$ such that
\begin{align*}
\int\left\vert \nabla v\right\vert ^{2}-\lambda\int gf^{2}(v) \ge\kappa \int\left\vert \nabla
v\right\vert ^{2}\text{,}%
\end{align*}
therefore the additional term $\int gf^{2}(v)$ in the functional does not
affect the verification of coerciveness. Moreover, similar to the functional
$\psi$ defined in (\ref{kq}), the functional $v\mapsto\int
gv^{2}$ is also weakly continuous on $D^{1,2}(\mathbb{R}^{N})$; therefore (\ref{e}) remains
valid even there is an additional term involving $g$ in the argument.

Similarly, because of the continuous embedding $H^{1}(\mathbb{R}^{N})\hookrightarrow D^{1,2}(\mathbb{R}^{N})$,
replacing the space $D^{1,2}(\mathbb{R}^{N})$ by $H^{1}(\mathbb{R}^{N})$ in the argument, we can obtain
similar results for%
\[
\left\{
\begin{array}
[c]{l}%
-\Delta u+V(x)u-u\Delta(u^{2})=k(x)\left\vert u\right\vert ^{q-2}%
u-h(x)\left\vert u\right\vert ^{s-2}u\text{,}\\
u\in H^{1}(\mathbb{R}^{N})\text{,}%
\end{array}
\right.
\]
where $V$ is a positive potential bounded away from $0$.
\end{rem}

\end{document}